\documentclass[12pt]{amsart}
\usepackage{amscd,amssymb}
\usepackage{tikz}
\usepackage[plainpages,backref,urlcolor=blue]{hyperref}
\usepackage[all]{xy}

\theoremstyle{plain}
\newtheorem{thm}[subsection]{Theorem}
\newtheorem{lem}[subsection]{Lemma}
\newtheorem{prop}[subsection]{Proposition}
\newtheorem{cor}[subsection]{Corollary}

\theoremstyle{definition}
\newtheorem{rk}[subsection]{Remark}
\newtheorem{definition}[subsection]{Definition}
\newtheorem{ex}[subsection]{Example}

\numberwithin{equation}{section}
\newcommand{\OO}{{\mathcal O}}

\newcommand{\F}{{\mathcal F}}

\newcommand{\wJ}{\widehat{J}}

\newcommand{\al}{{\alpha}}

\newcommand{\C}{\mathbb{C}}
\newcommand{\PP}{\mathbb{P}}

\newcommand{\N}{\mathbb{N}}

\DeclareMathOperator{\indeg}{indeg}
\DeclareMathOperator{\Syz}{Syz}

\begin{document}

\title [Addition-deletion results for the minimal degree]
{Addition-deletion results for the minimal degree of a Jacobian syzygy of a union of two curves}

\author[Alexandru Dimca]{Alexandru Dimca$^{1}$}
\address{Universit\'e C\^ ote d'Azur, CNRS, LJAD, France and Simion Stoilow Institute of Mathematics,
P.O. Box 1-764, RO-014700 Bucharest, Romania}
\email{dimca@unice.fr}

\author[Giovanna Ilardi]{Giovanna Ilardi}
\address{Dipartimento Matematica Ed Applicazioni ``R. Caccioppoli''
Universit\`{a} Degli Studi Di Napoli ``Federico II'' Via Cintia -
Complesso Universitario Di Monte S. Angelo 80126 - Napoli - Italia}
\email{giovanna.ilardi@unina.it}

\author[Gabriel Sticlaru]{Gabriel Sticlaru}
\address{Faculty of Mathematics and Informatics,
Ovidius University
Bd. Mamaia 124, 900527 Constanta,
Romania}
\email{gabriel.sticlaru@gmail.com}

\thanks{\vskip0\baselineskip
\vskip-\baselineskip
\noindent $^1$This work has been partially supported by the Romanian Ministry of Research and Innovation, CNCS - UEFISCDI, grant PN-III-P4-ID-PCE-2020-0029, within PNCDI III}



\begin{abstract} 
Let $C:f=0$ be a reduced curve in the complex projective plane.
The minimal degree $mdr(f)$ of a Jacobian syzygy for $f$, which is the same as the minimal degree of a derivation killing $f$, is an important invariant of the curve $C$, for instance it can be used to determined whether $C$ is free or nearly free.
In this note we study the relations of this invariant $mdr(f)$ with a decomposition of $C$ as a union of two curves $C_1$ and $C_2$, without common irreducible components. When all the singularities that occur are quasihomogeneous, a result by Schenck, Terao and Yoshinaga yields finer information on this invariant in this setting.
Using this, we give some geometrical criteria, {\it the first ones of this type in the existing literature as far as we know}, for a line to be a jumping line for the rank 2 vector bundle of logarithmic vector fields along a reduced curve $C$.
\end{abstract}
 
\maketitle


\section{Introduction} 
Let $S=\C[x,y,z]$ be the polynomial ring in three variables $x,y,z$ with complex coefficients. We denote by $ \partial_x, \partial_y, \partial_z$ the partial derivations with respect to $x,y,z$ respectively.
Let $Der(S)=\{\partial =a \partial_x + b  \partial_y+ c \partial_z \ : \ a,b,c \in S \}$ be the free $S$-module of $\C$-derivations of the polynomial ring $S$.
\begin{definition}
\label{def2} Let $g \in S$ be a polynomial.
The $S$-module $D(g)$ of derivations of $S$ {\it preserving the principal ideal} $(g) \subset S$ is by definition
$$D(g)=\{\partial \in Der(S) \ : \ \partial g \in (g) \}.$$
Moreover, the $S$-module $D_0(g)$ of derivations of $S$ {\it  killing the polynomial} $g$  is by definition
$$D_0(g)=\{\partial \in Der(R) \ : \ \partial g =0 \}.$$
\end{definition}
When $g$ is a homogeneous polynomial, then both modules $D(g)$ and $D_0(g)$ are graded $S$-modules and one has
$$D(g)=D_0(g) \oplus S(-1) \cdot E,$$
 where $E=x\partial_x+y\partial_y+z\partial_z$ denotes the {\it Euler derivation.} The curve $C_g: g=0$ in $\PP^2$ is said to be {\it  free} if 
 $D(g)$ or, equivalently, $D_0(g)$ is a free graded $S$-module.
 
 Note that $D_0(g)$ can be identified with the $S$-module of all Jacobian relations for $g$, namely to
$$AR(g)=\{(a,b,c) \in S^3 \ : \ ag_x+bg_y+cg_z=0 \},$$
where $g_u=\partial_u g$, for $u=x,y,z$. An important numerical invariant associated to a reduced curve $C:f=0$ in the projective plane $\PP^2$ is the {\it minimal degree of a derivation killing $f$} or,
equivalently, the {\it minimal degree of a Jacobian relation (syzygy) for $f$}. This is defined by
$$mdr(f)= \min \{s \in \N: \ D_0(f)_s  \ne 0 \}=\min \{s \in \N: \ AR(f)_s  \ne 0 \}.$$
It can be used for instance to characterize the free or the nearly free curves, see \eqref{F1} and \eqref{F2} below.
In this note we study the relations of this invariant with a decomposition of $C$ as a union of two curves $C_1:f_1=0$ and $C_2:f_2=0$, without common irreducible components. In particular, we would like to relate $r=mdr(f)=mdr(f_1f_2)$ to $r_j=mdr(f_j)$ for $j=1,2$. The case when $C_1$ is a  a line arrangement and $C_2$ is a line was studied in detail in \cite{ADS}.

In section 2 we recall some basic notations and facts, for instance the definition of the Jacobian module $N(f)$ and of free, nearly free and plus one generated curves which play a key role in this paper.

Then we consider in section 3 the case when $C_2$ is a line and $C_1$ is  any reduced curve, not having $C_2$ as a component. We study in Theorem \ref{thm1} the behavior of our invariant $mdr(f)$ when $C_2$ is a member of a pencil of lines in $\PP^2$, under the assumption that we know not only $r_1$, but also a non trivial derivation in $D_0(f_1)_{r_1}$. Several examples are given in section 4. 

We study in section 5 the general case of two curves $C_1$ and $C_2$, and get bounds for $r=mdr(f)$ in terms of the degrees $d_j=\deg(f_j)$  and of the invariants $r_j$ for $j=1,2$, see Theorem \ref{thm10}. As an example, we discuss in Proposition \ref{prop10} all the possibilities when both $C_1$ and $C_2$ are smooth conics.

Finally, in section 6, {\it we assume that all the singularities of $C_1$ and $C$ are quasihomogeneous and that $C_2$ is a smooth curve (most of the time $C_2$ is also supposed to be rational).}
Under this assumption, we may use a key result by Schenck, Terao and Yoshinaga, see \cite{STY}, to get finer information on $r$. Our Theorem \ref{thm20} gives a description of the cohomology exact sequence associated to the short sheaf exact sequence obtained in \cite{STY}, paying special attention to the description of the morphisms between the corresponding cohomology groups.

This approach was already used in  \cite{STY} to relate the freeness of $C_1$ to the freeness of $C$. Here we show that even when the curve $C_1$ is not free, one can obtain valuable information on $r$ using this approach.
This idea works best when the Jacobian module $N(f_1)$ is small, and this explains why we consider mostly free and nearly free curves $C_1$.
Sometimes the determination of $r$ is rather easy, using just the knowledge of the numerical invariant $r_1$, as in most examples in section 6.
 In  Example \ref{ex21E} we present a situation where one needs to use the morphisms in the exact sequence given by Theorem \ref{thm20}, namely the multiplication by $f_2^2$ between the two Jacobian modules $N(f_1)$ and $N(f)$.
As a by-product, under the assumption for this final section,
we get lower bounds on the number of points in the intersection $C_1 \cap C_2$ in terms of $r_1$ when the curve $C_1$ is  free or nearly free and $C_2$ is either a line or a smooth conic, see Corollary \ref{corF} and  Corollary \ref{cor21NF}.

Finally we give some geometrical criteria, {\it the first ones of this type in the existing literature as far as we know}, for a line $L$  in $\PP^2$ to be a jumping line for the rank 2 vector bundle  $T\langle C \rangle $  of logarithmic vector fields along a reduced curve $C$, see Theorem \ref{thmJ} and Example \ref{ex1.1} where this is applied to Thom-Sebastiani curves.

\medskip

We would like to thank Laurent Bus\'e, Piotr Pokora, \c Stefan Toh\u aneanu and Masahiko Yoshinaga for useful discussions related to this paper.

\section{Prerequisites}

In this section we recall some basic facts, see for instance  \cite{CDI,DSt3syz}. For any degree $e$ reduced homogeneous polynomial $g \in S_e$,  let $N(g)= \wJ_g/J_g$ be the Jacobian module of $g$, with $J_g$ the Jacobian ideal of $g$ in $S$, spanned by the partial derivatives $g_x,g_y$, $g_z$ of $g$, and $\wJ_g$ the saturation of the ideal $J_g$ with respect to the maximal ideal 
${\bf m}=(x,y,z)$ in $S$. 
We set   $n(g)_j=\dim N(g)_j$,  $T_g=3(e-2)$ and recall that  we have 
\begin{equation}\label{LPN} 
n(g)_0\leq n(g)_1\leq\ldots\leq n(g)_{\lfloor\frac{T_g}{2}\rfloor-1}\leq n(g)_{\lfloor\frac{T_g}{2}\rfloor}\geq n(g)_{\lfloor\frac{T_g}{2}\rfloor+1}\geq\ldots\geq n(g)_T.\end{equation}
For a reduced curve $C_g:g=0$, we consider the following  invariants   $$\sigma(C_g)=\min\left\{j : n(g)_j\neq 0\right\}=\indeg(N(f)) 
\text{ and }
\nu(C_g)=\max\left\{n(g)_j\right\}_j.$$ 
The self duality of the graded $S-$module $N(g)$ implies 
$
 n(g)_j=n(g)_{T_g-j}, 
$  
for any integer $j$,  see \cite{Se}. In particular $n(g)_k>0$ exactly when $\sigma(C_g) \leq k \leq T_g-\sigma(C_g).$\\ 
The  form of the minimal graded free resolution for the Milnor algebra $M(g)=S/J_g$ 
is 
\begin{equation}\label{RMM}0\to \displaystyle{\oplus_{i=1}^{m-2}S(-e_i)}\to \displaystyle{\oplus_{i=1}^{m}S(1-e-d'_i)}\to S^3(1-e)\to S,\end{equation} 
with $e_1\leq e_2\leq\ldots\leq e_{m-2}$ and $1 \leq d'_1\leq d'_2\leq\cdots\leq d'_m$. In this case the curve $C_g$ is said to be an $m$-{\it syzygy curve}. The first degree $r_g=d_1'$ is denoted by $mdr(g)$ and is the {\it minimal degree of a Jacobian relation (syzygy) for $g$}.
It follows from \cite[Lemma 1.1]{HS12} that one has $$e_j=e+d'_{j+2}-1+\epsilon_j,$$ for $j=1,\ldots, m-2$ and some integers $\epsilon_j\geq 1.$
The minimal resolution of $N(g)$ obtained from (\ref{RMM}), by \cite[Proposition 1.3]{HS12}, is
 $$0\to \displaystyle{\oplus_{i=1}^{m-2}S(-e_i)}\to\displaystyle{\oplus_{i=1}^mS(-\ell_i)}\to\displaystyle{\oplus_{i=1}^m}S(d'_i-2(e-1))\to \displaystyle{\oplus_{i=1}^{m-2}S(e_i-3(e-1))},$$
 where $\ell_i=e+d'_i-1$. It follows that 
\begin{equation}\label{sigma} 
 \sigma(C_g)=3(e-1)-e_{m-2}=2(e-1)-d'_m-\epsilon_{m-2}.
\end{equation} 
The following are important special cases, see \cite{Abe18,DStRIMS,DSt3syz}.
Here $\tau(C_g)$ is the total Tjurina number of the curve $C_g$, which is the same as the degree of the Jacobian ideal $J_g$.
\begin{enumerate}

\item $C_g$ is a free curve if and only if $m=2$ and $d_1'+d'_2=e-1$. In this case $\nu(C_g)=0$ and  $N(g)=0$. The degrees $(d_1',d_2')$ are the {\it exponents} of the free curve $C_g$.
Moreover, a reduced curve $C_g$ is free if and only if 
\begin{equation}\label{F1} 
 \tau(C_g)=(e-1)^2-r_g(e-r_g-1),
\end{equation} 
see \cite{Dmax,dPW}. 

\item $C_g$ is a nearly free curve if and only if $m=3$ and $d_1'+d'_2=e$, $d_3'=d_2'$. In this case $\nu(C_g)=1$ and  $\sigma(C_g)=e+d_1'-3$. The degrees $(d_1',d_2')$ are the {\it exponents} of the nearly free curve $C_g$.
Moreover, $C_g$ is nearly free if and only if 
\begin{equation}\label{F2} 
 \tau(C_g)=(e-1)^2-r_g(e-r_g-1)-1,
\end{equation} 
see \cite{Dmax}.

\item $C_g$ is a plus one generated curve if and only if $m=3$ and $d_1'+d'_2=e$, $d_3'>d_2'$, see \cite{Abe18} for the case $C_g$ a line arrangement and \cite{DSt3syz} for the general case. In this case $\nu(C_g)=d_3'-d_2'+1$ and   $\sigma(C_g)=2e-d_3'-3$.

\end{enumerate}

\section{Adding a line to a reduced curve}

Consider a reduced plane curve $C_1:f_1=0$  of degree $d_1$ in $\PP^2$ such that $mdr(f_1)=r_1$. Let $L$ be a line in $\PP^2$, which is not an irreducible component of $C_1$ and consider the curve $C=C_1 \cup L: f=0$.
Then $C$ has degree $d=d_1+1$, and we denote $r=mdr(f)$.
In this section we analyze the relation between $r$ and $r_1$, starting with the following result. 
\begin{prop}
\label{prop1}
With the above notation, one has $r_1 \leq r \leq r_1+1$.
\end{prop}
\proof
Choose a coordinate system on $\PP^2$ such that the line $L$ is given by $z=0$, and hence $f=zf_1$. Let
\begin{equation}
\label{syzf}
af_x+bf_y+cf_z=0
\end{equation}
be a Jacobian syzygy of minimal degree $r$ for $f$, and
\begin{equation}
\label{syzf'}
a_1f_{1x}+b_1f_{1y}+c_1f_{1z}=0
\end{equation}
 a Jacobian syzygy of minimal degree $r_1$ for $f_1$.
Note that one has
\begin{equation}
\label{e1}
f_x=zf_{1x}, \  f_y=zf_{1y} \text{ and } f_z=zf_{1z}+f_1=\frac{1}{d_1}xf_{1x}+\frac{1}{d_1}yf_{1y} +\frac{d}{d_1}zf_{1z}.
\end{equation}
Using \eqref{syzf} we get
$$azf_{1x}+bzf_{1y}+c(zf_{1z}+f_1)=0,$$
and hence the polynomial $c$ is divisible by $z$, so we can write $c=zc'$. Indeed, note that $f_1$ is not divisible by $z$ by our assumptions. With this notation, and using \eqref{e1}, we get after division by $z$ the following equation.
\begin{equation}
\label{e2}
(a+\frac{1}{d_1}c'x)f_{1x}+(b+\frac{1}{d_1}c'y)f_{1y}+\frac{d}{d_1}c'zf_{1z}=0.
\end{equation}
This implies $r_1 \leq r$.
Similarly, using \eqref{syzf'} and \eqref{e1} we get
\begin{equation}
\label{e3}
(d_1a_1z-c_1x)f_x+(d_1b_1z-c_1y)f_y+c_1d_1zf_z=0.
\end{equation}
Note that this is a non trivial syzygy, namely one cannot have
$$d_1a_1z-c_1x= d_1b_1z-c_1y=c_1d_1z= 0.$$
This implies  $r \leq r_1+1$.
\endproof
\begin{rk}
\label{rk01}
With the above notation, if $z$ divides  $c_1$, the coefficient of $f_{1z}$ in \eqref{syzf'}, then all the coefficients in \eqref{e3} are divisible by $z$,
and hence after simplification by $z$ we get $r=r_1$ in this case. When $\dim D_0(f_1)_{r_1}>1$, there is a choice of the syzygy \eqref{syzf'} within a linear system, and some choices may be better than others, i.e. for the good ones $z$ divides  $c_1$, see Example \ref{exPOG} below for such a situation.
\end{rk}
To say more about
the value of $r$, it is convenient to look not only at a single line $L$, but at all the lines in a pencil. 
The pencil we consider is formed by all the lines in $\PP^2$ passing through a point $p \in \PP^2$, which may or may not be on the curve $C_1$. We choose a coordinate system on $\PP^2$ such that $p=(1:0:0)$,
hence a line in the pencil has the equation $L_u: sy+tz=0$ for some
$u=(s:t) \in \PP^1$. Assume that \eqref{syzf} and 
\eqref{syzf'} are minimal degree Jacobian syzygies for $f=(sy+tz)f_1$ and
respectively for $f_1$, with respect to this coordinate system. 
Note that the coefficients $a_1,b_1,c_1$ are known and independent of $u$, since they depend only on $C_1$ and the choice of the coordinate system.
Let 
$$r=d'_1(f) \leq d'_2(f) \leq \dots \leq d'_m(f)$$
 be the degrees of a minimal set of generators for $AR(f)$ coming from the resolution \eqref{RMM} of the Milnor algebra $M(f)$, which depend in general on $u$, see Example \ref{ex2} below. Elementary computations similar to those done above yield the following syzygy
\begin{equation}
\label{e4}
A_uf_x+B_uf_y+C_uf_z=0,
\end{equation}
where $A_u=d(sy+tz)a_1-x(sb_1+tc_1)$, $B_u=d(sy+tz)b_1-y(sb_1+tc_1)$ and finally
$C_u=d(sy+tz)c_1-z(sb_1+tc_1)$.
Using this syzygy, we can prove the following result.
\begin{thm}
\label{thm1}
With the above notation, if $sy+tz$ is a factor of $sb_1+tc_1$, then $r=r_1$.
If $sy+tz$ is not a factor of $sb_1+tc_1$, then either
\begin{enumerate}
\item $r=r_1+1$, or 

\item $r=r_1$ and $d'_2(f) \leq r+1$.
\end{enumerate}
Moreover, the case $(2)$ is impossible when $2r_1<d_1-1$, or when $2r_1=d_1-1$ and $C$ is not free.

\end{thm}

\proof
The first claim is obvious. Indeed, when $sy+tz$ is a factor of $sb_1+tc_1$, the coefficients $A_u$, $B_u$ and $C_u$ can be divided by $sy+tz$, and the syzygy \eqref{e4} yields a syzygy of degree $r_1$. Since $r \geq r_1$ by Proposition \ref{prop1}, we get $r=r_1$.
Assume now that $sy+tz$ is not a factor of $sb_1+tc_1$. Then we claim that
the syzygy \eqref{e4} is primitive, i.e. it is not a multiple of a syzygy of strictly lower degree. In other words, we have to show that $A_u$, $B_u$ and $C_u$ have no common factor in this case.
Note that $yA_u-xB_u=d(sy+tz)(ya_1-xb_1)$, $zA_u-xC_u=d(sy+tz)(za_1-xc_1)$ and
$zB_u-yC_u=d(sy+tz)(zb_1-yc'_1$.
Let $D$ be a common irreducible factor of $A_u$, $B_u$ and $C_u$, supposed to be a homogeneous polynomial of degree $>0$. It it clear that $D$ cannot be $sy+tz$, since $sy+tz$ is not a factor of $sb_1+tc_1$. Hence $D$ has to divide the polynomials $m_{12}=ya_1-xb_1$, $m_{13}=za_1-xc_1$ and
$m_{23}=zb_1-yc_1$.
Recall now the construction of the Bourbaki ideal $B(C_1, \rho'_1)$ associated to the curve $C_1$ and to the minimal degree syzygy $\rho'_1$
given by \eqref{syzf'}, as described in \cite[Section 5]{DStJump}.
It follows that the Bourbaki ideal $B(C_1, \rho'_1)$ is contained in the principal ideal generated by $D$. This is a contradiction, since 
the Bourbaki ideal $B(C_1, \rho'_1)$ defines a subscheme which is either empty (when $C_1$ is a free curve), or zero-dimensional, see
\cite[Theorem 5.1]{DStJump}.

Therefore the syzygy \eqref{e4} is indeed primitive.
It follows that either $r=r_1+1$, or $r=r_1$ and $d_2(f) \leq r+1$.
Note that in this latter case we have
$$d_1=d-1 \leq d'_1(f) +  d'_2(f) \leq r_1+r_1+1=2r_1+1.$$
Indeed, recall that 
$d-1=d'_1(f)+d'_2(f)$
exactly when $C$ is free, and 
$d-1<d'_1(f)+d'_2(f)$ otherwise, see for instance \cite{ST}.
\endproof
\begin{prop}
\label{prop01} 
With the notation from Theorem \ref{thm1}, we have the following equivalent properties.
\begin{enumerate}
\item $sy+tz$ is a factor of $sb_1+tc_1$ for infinitely many $u=(s:t) \in \PP^1$;

\item $sy+tz$ is a factor of $sb_1+tc_1$ for all $u=(s:t) \in \PP^1$;

\item the reduced curve $C_1:f_1=0$ is the union of the curve $h=0$ with a pencil of lines $g=0$ passing through the point $p=(1:0:0)$.

\end{enumerate}

\end{prop}

\proof
The fact that (2) implies (1) is clear. First we show that (1) implies (3).
Note that $sy+tz$ is a factor of $sb_1+tc_1$ for infinitely many $u=(s:t) \in \PP^1$ if and only if there is a polynomial $h$ of degree $r_1-1$ such that
$b_1=yh$ and $c_1=zh$. Replacing these values in \eqref{syzf'} we conclude that $f_{1x}$ is divisible by $h$, say $f_{1x}=hg$, with $\deg g=d_1-r_1\geq 1$. If we divide the syzygy \eqref{syzf'} by $h$, we get
\begin{equation}
\label{e5}
a_1g+yf_{1y}+zf_{1z}=0
\end{equation}
or, equivalently,
\begin{equation}
\label{e6}
a_1g+d_1f_1-xf_{1x}=0.
\end{equation}
It follows that $g$ is a common factor of $f_1$ and $f_{1x}$. To conclude the proof of the implication $(1) \implies (3)$ we use the following result, communicated to us by Laurent Bus\' e.

\begin{lem}
\label{lempropBL}
With the above notation, assume that $g=G.C.D.(f_1,f_{1x})$ has degree $\geq 1$. Then $g$ is a homogeneous polynomial in $y$ and $z$ only, and
$$f_1(x,y,z)=g(y,z)h(x,y,z),$$
for some homogeneous polynomial $h \in S$. In geometric terms, the reduced curve $C_1:f_1=0$ is the union of the curve $h=0$ with a pencil of lines $g=0$ passing through the point $p=(1:0:0)$.
\end{lem}

\proof

Let $A$ be an irreducible common factor of $f_1$ and $f_{1x}$, such that
$f_1=A U$ for $U \in S$. This implies
$f_{1x}=A_xU+AU_x$, and hence, if $A_x \ne 0$, then $A$ has to divide
$U$. Indeed, $A$ cannot divide $A_x$ since $\deg A_x <\deg A$.
But this contradicts the fact that $C_1:f_1=0$ is a reduced curve.
Hence $A_x=0$, in other words $A$ is a homogeneous polynomial in $y$ and $z$ only. Since $g$ is a product of such polynomials, the claim is proved.
\endproof

Finally we show that (3) implies (2). Assume that $g=G.C.D.(f_1,f_{1x})$ has degree $\geq 1$, then
one can define $a_1$ using the above equation \eqref{e6}. Then, if we multiply the equation  \eqref{e5} by $h=f_{1x}g^{-1}$, we get a primitive syzygy of the form \eqref{syzf'}, where $b_1=yh$ and $c_1=zh$. 
\endproof

\section{Examples}

\begin{ex}
\label{ex1}
Assume $C_1$ is an irreducible nodal curve and $L_u:sy+tz=0$ is a line such that $C=C_1 \cup L_u$ is nodal.
Then it is known that $r_1=d_1-1$ and $r=d-2=d_1-1$, see \cite{DStEdin, E}. Note that one has in this case $d_2(f)=r+1$,  see \cite[Theorem 4.1] {DStEdin}.
Hence  the case (2) of  Theorem \ref{thm1} might occur.
\end{ex}

\begin{ex}
\label{ex2}
Consider the rational cuspidal curve $C_1:f_1=xy^{d_1-1}+z^{d_1}=0$, $d_1 \geq3$, which is nearly free, and $L_u:sy+tz=0$ a line passing through the singular point $p=(1:0:0)$.
Then the syzygy \eqref{syzf'} becomes
$$(d_1-1)xf_{1x}-yf_{1y}=0.$$
Hence $sb_1+tc_1=-sy$ is divisible by $sy+tz$ only for $(s:t)=(1:0)$ and for $(s:t)=(0:1)$, and we get in these cases $r=r_1=1$ using Theorem \ref{thm1} as we see now.
The curve $C': f=xy^{d_1}+yz^{d_1}=0$ corresponding to $(s:t)=(1:0)$
is free, the two generating syzygies being
$$(d_1)^2xf_x-d_1yf_y+zf_z=0$$
and
$$d_1z^{d_1-1}f_x-y^{d_1-1}f_z=0$$
satisfying $d'_1(f)+d'_2(f)=1+ (d_1-1)=d-1$. 
The curve $C'': f=xy^{d_1-1}z+z^{d_1+1}=0$ corresponding to $(s:t)=(0:1)$,
is nearly free with exponents $d'_1(f)=1$, $d'_2(f)=d'_3(f)=d-1$. Indeed, note that the curve $C''$ has two singularities, namely $p=(1:0:0)$ and
$q=(0:1:0)$. The singularity at $q$ is a simple node $A_1$, and the singularity at $p$ is given in local coordinates $y'=y/x$ and $z'=z/x$ by
$(y')^{d_1-1}z'+(z')^{d_1+1}=0$. This is a quasi homogeneous singularity,
with weights $wt(z')=d^{-1}$ and $wt(y')=d_1[(d_1-1)d]^{-1}$.
It follows that 
$$\tau(C'',p)=\mu(C'',p)=d^2-3d+1$$
and hence the total Tjurina number of $C''$ is given by
$$\tau(C'')=\tau(C'',q)+ \tau(C'',p)=d^2-3d+2.$$
The fact that $C''$ is nearly free follows now from \eqref{F2}.

For $d_1 \geq 4$ and for $L_u:y+z=0$, we have $r=r'+1=2$
by Theorem \ref{thm1}, since $2r_1<d_1-1$ in this case. 
The corresponding curves $C_u$ are again nearly free, but this time with exponents 
$d'_1(f)=2$, $d'_2(f)=d'_3(f)=d-2$. To see this, one notes that a curve $C_u$ in this family has two singularities, a node and a semi quasi homogeneous singularity  $(C_u,p):g(y',z')=g_0(y',z')+g_+(y',z')=0$,
where $g_0$ is quasi homogeneous and $g_+(y',z')$ is the sum of two monomials of strictly higher degree. Working in the Milnor algebra $M(g_0)$, we see that the Tjurina algebra of $g$ is isomorphic to the quotient $M(g_0)/(y^{d_1})$. This implies that 
$$\mu(C_u,p)=(d_1)^2-d_1-1 \text{ and } \tau(C_u,p)=(d_1-1)^2+1.$$
It follows that $\tau(C_u)=(d_1-1)^2+2=(d-2)^2+2=(d-1)^2-2(d-3)-1$,
showing that $C_u$ is nearly free by \eqref{F2}.

\end{ex}

\begin{ex}
\label{ex3}
Let $C_1:f_1=(y^2-2xy+z^2)(y^2+4xy+z^2)=0$, be the union of two smooth conics tangent at one point $p=(1:0:0)$ and meeting transversely at $q_{\pm}=(0:1:\pm i)$. Then using Singular we see that $r_1=2$ and a minimal degree derivation is given by
$$\partial'=xz\partial_x-yz\partial_y+y^2\partial_z.$$
Then the equation  $sy+tz$ of a line $L$ passing through the tangency point $p$ divides 
$$sb_1+tc_1=y(ty-sz)$$
if and only if either $(s:t)=(1:0)$ or $(s:t)=(1:\pm i)$.

The case $(s:t)=(1:0)$ corresponds to a common tangent $y=0$ to the two conics at $p$. Using Singular, we see that the corresponding curve 
$$C: f=yf_1=y(y^2-2xy+z^2)(y^2+4xy+z^2)=0$$
is free with exponents $(2,2)$, in particular $r=2=r_1$ as predicted by Theorem \ref{thm1}. 

The case $(s:t)=(1:\pm i)$ corresponds to a line joining the tangency point $p$ to one of the two nodes $q_{\pm}$ of $C_1$.
Using Singular, we see that the corresponding curve 
$$C: f=(y \pm iz)f_1=(y \pm iz)(y^2-2xy+z^2)(y^2+4xy+z^2)=0$$
is nearly free with exponents $(2,3)$, in particular, again $r=2=r_1$ as predicted by Theorem \ref{thm1}.

Finally, to see what happens when $sy+tz$ does not  divide 
$sb_1+tc_1=y(ty-sz)$, namely when the line through $p$ is general, we consider the special case $(s:t)=(1:1)$. Using Singular, we see that the corresponding curve 
$$C: f=(y +z)f_1=(y +z)(y^2-2xy+z^2)(y^2+4xy+z^2)=0$$
is a maximal Tjurina curve of type $(d,r)=(5,3)$, see \cite{DStmax} for the definition and the properties of such curves, and in particular $C$ has exponents  $(3,3,3,3)$. Hence $r=3=r_1+1$. Note that we can show that for any line  $L:sy+tz=0$ with $t \ne 0$, the singularity of $C$ at $p$ is of type $D_6$. Indeed, it follows easily that this singularity is semi weighted homogeneous of type
$(2,1;5)$, where $wt(y)=2$ and $wt(z)=1$. The claim follows using \cite[Corollary (7.39)]{RCS}. In particular, when $sy+tz$ does not  divide 
$sb_1+tc_1=y(ty-sz)$, we always have $\tau(C)=10$, since there are 4 nodes $A_1$  on $C$ in addition to the $D_6$ singularity.

\end{ex}

\begin{ex}
\label{exPOG}
Let $C_1: f_1=(y^2-xz)^2+y^2z^2+z^4=0$ be the curve considered in 
 \cite[Example 4.1]{DSt3syz}. This curve is plus one generated with exponents $(d_1',d_2',d_3')=(2,2,3)$, in particular $\dim D_0(f_1)_2=2$.
 If we choose the right element in $D_0(f_1)_2$, namely
 $$ \partial'=(2xy+3yz)\partial_x+(xz+2z^2)\partial_y-yz\partial_z,$$
 then $z$ divides the coefficient of $\partial_z$, and 
it follows that $r=r_1=2$ by Theorem \ref{thm1}.

\end{ex}

\section{The general case: the union of two curves}

Let $C_1:f_1=0$ and $C_2:f_2=0$ be two reduced curves in $\PP^2$, without common irreducible components. We denote $d_j=\deg f_j$ and $r_j=mdr(f_j)$ for $j=1,2$. Consider now the union of the two curves
$C:f=f_1f_2=0$, and let $d=d_1+d_2=\deg f$ and $r=mdr(f)$.
\begin{thm}
\label{thm10}
With the above notation, one has the following.
\begin{enumerate}
\item If $\delta_1 \in D_0(f_1)$, then
$$\delta=f_2\delta_1-\frac{\delta_1(f_2)}{d}E \in D_0(f),$$
where $E=x\partial_x+y\partial_y+z\partial_z$ denotes the Euler derivation. In particular
$$r\leq \min\{r_1+d_2, r_2+d_1\}.$$

\item  $D_0(f) \subset D(f_1)\cap D(f_2)$. More precisely,  for $\delta \ne 0$, one has $\delta \in D_0(f)$  if and only if $\delta$ can be written in a unique way in the form
$$ \delta=\frac{h}{d_1}E+ \delta_1 = -\frac{h}{d_2}E+ \delta_2, $$
where $h \in S$ and $\delta_j \in D_0(f_j)$ are non-zero derivations. 
In particular
$$r \geq \max \{r_1, r_2\}.$$
\end{enumerate}

\end{thm}
\proof
To prove $(1)$, first we check that $\delta( f)=0$. Then we note that 
$\delta \ne 0$ if $\delta_1 \ne 0$. Indeed, if $\delta_1(f_2)=0$, then clearly $\delta =f_2\delta_1\ne 0$. When $\delta_1(f_2)\ne 0$, note that
$$\delta (f_1)=\frac{d_1f_1\delta_1(f_2)}{d} \ne 0.$$
The last claim follows by noting that if $\delta_1$ is a homogeneous derivation then also $\delta$ is a homogeneous derivation. Moreover, the roles played by $f_1$ and $f_2$ are symmetric.

To prove $(2)$, 
start with $\delta \in D_0(f)$ and hence
$$\delta (f)=f_2\delta (f_1)+f_1\delta (f_2)=0.$$
If $\delta (f_1)=0$, then $\delta (f_2)=0$ and hence $\delta \in D_0(f_1)\cap D_0(f_2)$.
If $\delta( f_1)\ne 0$, then $f_2$ divides the product $f_1\delta (f_2)$. Since
$f_1$ and $f_2$ have no common factor by our assumptions, it follows that $f_2$ divides $\delta (f_2)$, hence $\delta \in D(f_2)$. This is possible only if   $\delta \in D(f_1)$ as well. It follows that we  can write $\delta \in D_0(f)$  in the form
$$\delta=h_jE+\delta_j$$ 
where $h_j \in S$ and $\delta_j \in D_0(f_j)$. Clearly $\delta_j \ne 0$, since otherwise $\delta (f) \ne 0$. 
Then $\delta (f_1)= d_1h_1f_1$ and $\delta (f_2)=d_2h_2f_2$.
It follows that
$$0=\delta (f)=\delta (f_1)f_2+f_1\delta (f_2)=f_1f_2(d_1h_1+d_2h_2).$$
Then one implication in the claim follows by taking $h=d_1h_1=-d_2h_2$. The other implication is obvious.
\endproof

\begin{rk}
\label{rk10}The inequality $r \geq \max \{r_1, r_2\}$ was already noticed in \cite[Proposition 3.2. (ii)]{BST}, where the $S$-module $D_0(g)=AR(g)$ is denoted by $\Syz (J_g)$ and $mdr(g)$ is denoted by
$\indeg (\Syz (J_g))$. Note also that in \cite{BST} one works over the polynomial ring in $n$-variables with coefficients in an arbitrary infinite field.
The corresponding result for a product $f=f_1f_2 \cdot \ldots \cdot f_m$ of $m\geq 2$ forms in $n$-variables is considered in \cite[Proposition 3.5]{BST}. Interesting information on the invariant  $\indeg (\Syz (J_f))$ when  $C:f=0$ is the union of several smooth plane curves meeting transversally is given in \cite[Proposition 3.6]{T}.

\end{rk}

\begin{cor}
\label{cor10}
With the above notations, $r=mdr(f)$ is the minimal integer $s$ such that
either $D_0(f_1)_s\cap D_0(f_2)_s \ne 0$, or $D_0(f_1)_s+D_0(f_2)_s$ contains a non-zero multiple of the Euler derivation $E$.
\end{cor}
\proof
The first case corresponds to $h=0$ in Theorem \ref{thm10} (2), while the second case corresponds to $h \ne 0$.
\endproof

\begin{ex}
\label{ex10}
Let $C_1:f_1=x^2 + y^2 -z^2=0$ and $C_2:f_1=x^2 + y^2 -4z^2=0$ be  two smooth conics with 2 tacnodes as in Proposition \ref{prop10} $(3)$. Hence $d_1=d_2=2$, $r_1=r_2=1$.
Note that  $y\partial_x-x\partial_y \in D_0(f_1)_1 \cap D_0(f_2)_1$. Therefore, according to Corollary \ref{cor10} we have $r=1$, see also Proposition \ref{prop10}, (3).\\
Consider next the case  $C_1:f_1=xyz=0$ and $C_2:f_2=xy+yz+xz=0$.
Then $C_2$ is a smooth conic circumscribed in the triangle $C_1$.
Using Singular, we see that $r=2$ and $D_0(f)_2$ is spanned by
$$\delta=2x(y-z)\partial _x-y(3y+2z)\partial_y+z(2y+3z)\partial_z$$
and
$$\delta'=x(3x+4y-2z)\partial_x-y(2x+6y+2z)\partial_y+z(-2x+4y+3z)\partial_z.$$
Then $\delta (f_1)= xyz(y+3z)=d_2h_2f_1$, which implies $h=-d_2h_2=-(y+3z)$.
Similarly $\delta' (f_1)= xyz(-x+2y-z)=d_2h_2f_1$, which implies that in this case $h=-d_2h_2=x-2y+z$. It follows that in this case
$D_0(f_1)_2 \cap D_0(f_2)_2=0$. Therefore, both situations may occur in Corollary \ref{cor10}. The fact that $r=2$ in this case is discussed from another view-point, without the use of Singular, in Example \ref{ex10bis}.
The curve $C$ has three $D_4$ singularities and hence $\tau(C)=12$.
Using the characterization of free curves in \eqref{F1}, it follows that $C$ is a free curve.
\end{ex}

Let $C_1$ and $C_2$ be smooth conics, hence $d_1=d_2=2$ and $r_1=r_2=1$. For $C=C_1 \cup C_2$, Theorem \ref{thm10} gives us $1 \leq r \leq 3$. We have the following precise result.

\begin{prop}
\label{prop10}
The two conics $C_1$ and $C_2$  can be in one of the following four situations. 
\begin{enumerate}

\item $|C_1 \cap C_2|=4$, and then all the intersection points are nodes for $C$.
 In this case  $r=2$.

\item  $|C_1 \cap C_2|=3$, and then one intersection point is a tacnode  and the other two intersection points are nodes for $C$. In this case $r=2$.

\item  $|C_1 \cap C_2|=2$. Then the two intersection points are either two tacnodes for $C$,  and in this case $r=1$ and the curve $C$ is nearly free with exponents $(1,3)$, or a node $A_1$ and a singularity $A_5$ for $C$,
and in this case $r=2$ and the curve $C$ is nearly free with exponents $(2,2)$.

\item  $|C_1 \cap C_2|=1$, and then the intersection point is a singularity $A_7$ for $C$,  $r=1$ and $C$ is a free curve.

\end{enumerate}
\end{prop}

Computations with Singular suggest that in case $(1)$ the curve  $C=C_1 \cup C_2$ is a 4-syzygy curve with
exponents $(2,3,3,3)$, and in case $(2)$ the curve $C=C_1 \cup C_2$ is a plus one generated curve with
exponents $(2,2,3)$.

\proof
The claim $(1)$ follows from \cite[Theorem 4.1]{DStEdin}.
For the claim $(2)$ we use the inequalities involving $r$ and the Tjurina number $\tau(C)$ due to du Plessis and Wall, see \cite{dPW}.
In case $(2)$ we have $\tau(C)=5$. We know that
$$5=\tau(C) \geq (d-1)(d-1-r)=3(3-r).$$
This implies $r \geq 2$. For $r=3$ we also have
$$5=\tau(C) \leq (d-1)(d-1-r)+r^2-{2r-d+2 \choose 2}= 3,$$
a contradiction. So the only possibility is $r=2$.
In case $(3)$, when the contact between $C_1$ and $C_2$ consists of two tacnodes, using the results in  \cite{Meg}, we see that a pair of conics in this situation is projectively equivalent to a pair of conics of the form
\begin{equation} \label{M1} 
C_1:f_1=x^2 - y^2 -z^2=0 \text{ and } C_2:f_1=x^2 - y^2 -k z^2=0,
\end{equation} 
with $k \in \C^*$, $k \ne 0$.
 We have $\tau(C)=6$, and the same approach as in Example \ref{ex10} gives $r=1$ in this case. Since we have $\tau(C)=6$ in this case, the equality
$$\tau(C) = (d-1)(d-1-r)+r^2-1,$$
holds, and
it follows from \eqref{F2} that $C$ is nearly free.
Assume now that the smooth conics $C_1$ and $C_2$ have a contact of type $A_5$  at $(0:0:1)$. Choosing the coordinates, we may assume that $C_1: f_1=yz-x^2=0$. Then it is easy to see that the other conic has an equation of the form $C_2=f_2=yz-x^2+a\cdot xy+b \cdot y^2=0$ for some $a \in \C^*$ and $b \in \C$. Since $\tau(C) =6$, we get as above $r \leq 2$. To show that $r>1$ one can use Theorem \ref{thm10} (2), since we have a simple description of $D(f_1)_1$. We see that $\delta_1 (f_2) \notin \C \cdot f_2$, for any $\delta_1 \in D_0(f_1)_1$. It follows that $r=2$ and we get the nearly freeness of $C$ as above. 

In the case (4), it follows from \cite[Proposition 1.3]{Bruce}, that the equations of the two conics can be chosen as follows
$$C_1: f_1=x(x+y)+yz-a \cdot y^2=0 \text{ and } C_2: f_2=x(x+y)+yz+a \cdot y^2=0$$
for some $a \in \C^*$. It follows that $(f_1)_x=(f_2)_x$, $(f_1)_z=(f_2)_z$, $\delta=(f_1)_z\partial_x-(f_1)_x\partial_z \in
D_0(f_1)_1 \cap D_0(f_1)_1 \ne 0$, and hence $r=1$. Moreover, the freeness of $C$ follows from
\cite{Dmax,dPW}.
\endproof

\begin{ex}
\label{ex11}
Let $C_1:f_1=(x^2 + y^2 -z^2)(x^2+y^2-4z^2)=0$ be the union of two smooth conics with 2 tacnodes as in Proposition \ref{prop10} $(3)$. Hence $d_1=4$, $r_1=1$.
Let $C_2:f_2= (x-z)(3y^2-(x+2z)^2)=0$ be the union of 3 lines, forming a triangle which is inscribed in the conic $C_2$ and circumscribed to the conic $C_1$. Then $d_2=3$, $r_2=1$.  For $C=C_1 \cup C_2$, Theorem \ref{thm10} gives us $1 \leq r \leq 3$.
Using Singular we see that $r=3$, see also Example \ref{ex11bis} below for a different approach. In fact, $C$ is a free curve with exponents $(3,3)$ as follows from \eqref{F1}, see also \cite{DPo}.
\end{ex}

\section{The case of quasihomogeneous singularities}

Consider the sheafification 
$$E_C:=
\widetilde{AR(f)}= \widetilde{D_0(f)} $$
of the graded $S$-module $AR(f)=D_0(f)$, which is a rank two vector bundle on $\PP^2$, see \cite{Se} for details. Moreover, recall that
\begin{equation} \label{equa1} 
E_C=T\langle C \rangle (-1),
\end{equation}
where $T\langle C \rangle $ is the sheaf of logarithmic vector fields along $C$ as considered for instance in \cite{FV1,MaVa,DS14}.
One has,  for any integer $k$, 
\begin{equation}
\label{e7}
H^0(\PP^2, E_C(k))=D_0(f)_k \text{ and }  H^1(\PP^2, E_C(k))=N(f)_{k+d-1},  
\end{equation}
where $d=\deg(f)$, for which we refer 
to \cite[Proposition 2.1]{Se}.
Return now to the setting of the previous section, where $C=C_1 \cup C_2$ and $f=f_1f_2$, and recall the following  result, see \cite[Theorem 1.6 and Remark 1.8]{STY}.

\begin{thm}
\label{thm20STY}
With the above notation, assume that $C_2$ is an irreducible curve, and that all singularities of $C_1$, $C_2$ and $C$ are quasihomogeneous. 
If $C_1 \cap C_2$ is contained in the smooth part of $C_2$, then there is an exact  sequence of sheaves on $ \PP^2$ given by
$$ 0 \to E_{C_1}(1-d_2) \stackrel{f_2} \longrightarrow  E_C(1) \to i_{2*}\F \to 0$$
 where  $i_2: C_2  \to \PP^2$ is the inclusion and $\F$ a torsion free sheaf  on ${C_2}$.
Moreover,
when $C_2$ is smooth, then one has $\F= \OO_{C_2}(-K_{C_2}-R)$, where $K_{C_2}$ is the canonical divisor on ${C_2}$ and $R$ is the reduced scheme of $C_1 \cap C_2$.
\end{thm}

{\it For simplicity, in this note we consider only the case $C_2$ smooth.}
If we set
$$\OO_{C_2}(1)=i_2^*\OO_{\PP^2}(1),$$
then one can write $\OO_{C_2}(1)=\OO_{C_2}(D)$, where the divisor $D$ corresponds to the intersection of a line in $\PP^2$ with the curve $C_2$, and hence $\deg D=d_2$. With this notation, by tensoring the above exact sequence with $\OO_{\PP^2}(k-1)$, for any integer $k$, we get the exact sequence
\begin{equation}
\label{e8}
0 \to E_{C_1}(k-d_2) \stackrel{f_2} \longrightarrow  E_C(k) \to i_{2*}\OO_{C_2}(-K_{C_2}-R+(k-1)D) \to 0.
\end{equation}
By taking the corresponding long cohomology sequence and using \eqref{e7}, we get the following result.
\begin{thm}
\label{thm20}
With the above notation, assume that $C_2$ is a smooth curve, and that all singularities of $C_1$ and $C$ are quasihomogeneous. Then there is an exact sequence  for any integer $k$ given by
$$ 0 \to D_0(f_1)_{k-d_2} \stackrel{\phi_k} \longrightarrow  D_0(f)_k \to H^0(C_2,\OO_{C_2}(-K_{C_2}-R+(k-1)D)) \to $$
$$ \to N(f_1)_{k-d_2+d_1-1} \stackrel{\psi_k} \longrightarrow N(f)_{k+d-1}   \to H^1(C_2,\OO_{C_2}(-K_{C_2}-R+(k-1)D)), $$
where the morphism $\phi_k: D_0(f_1)_{k-d_2} \to  D_0(f)_k$ is given by 
$$\phi_k(\delta_1)=f_2\delta_1-\frac{\delta_1(f_2)}{d}E$$
for $\delta_1 \in D_0(f_1)$ and $\psi_k$ is induced by the multiplication by $f_2^2$.
 In particular, if 
$$(k+2)d_2 <d_2^2+ |R|,$$
then the morphism $\phi_k$ is an isomorphism and $\psi_k$ is a monomorphism.
\end{thm}
\proof
The exact sequence above is part of the long cohomology exact sequence associated to the exact sequence of sheaves \eqref{e8}. It remains to explain the claims about the morphisms $\phi_k$ and $\psi_k$.
Using the identification $D_0(g)=D(g)/SE$, valid for any homogeneous polynomial $g \in S$, it is shown in \cite{STY} that the morphism $E_{C_1}(1-d_2) \to E_C(1)$ is induced by the multiplication by $f_2$.
In terms of the modules $D_0(g)$, this is precisely the mapping
$D_0(f_1)(-d_2) \to D_0(f)$ given by
$$\phi: \delta_1 \mapsto \delta=f_2\delta_1-\frac{\delta_1(f_2)}{d}E,$$
as constructed in Theorem \ref{thm10}. To explain the formula for $\psi_k$, consider the diagram of graded $S$-modules
$$
\xymatrix{0 \ar[r] & D_0(f_1)(-d_2) \ar[d]_{\phi}  \ar[r]^{\iota} & S^3(-d_2) \ar[d]_{\phi}  \ar[r]^{\nabla f_1} &
                            J_{f_1}(d_1-d_2-1) \ar[d]_{f_2^2}  \ar[r] & 0 \\
 0 \ar[r] & D_0(f)   \ar[r]^{\iota}  & S^3  \ar[r]^{\nabla f}  &    J_f(d-1)   \ar[r] & 0
}$$
Here $\iota$ are the obvious inclusions, $\phi$ is the morphism defined above and its extension to a map $Der(S)(-d_2)=S^3(-d_2) \to Der(S)=S^3$ given by the same formula, 
$\nabla f :S^3 \to J_f$ is the map $(a,b,c) \mapsto af_x+bf_y+cf_z$ and similarly for $\nabla f_1$, while 
$$f_2^2: J_{f_1}(d_1-d_2-1) \to J_f(d-1)$$
is the multiplication by $f_2^2$. A simple computation shows that this diagram is commutative. Since $N(f_1)= \wJ_{f_1}/J_{f_1}$ and $N(f)= \wJ_f/J_f$, it follows that the morphism $N(f_1)_{k-d_2+d_1-1} \stackrel{\psi_k} \longrightarrow N(f)_{k+d-1} $ induced by the long cohomology exact sequence, and hence coming from $\phi$, is nothing else but multiplication by $f_2^2$.  
The final inequality says that
$$\deg( -K_{C_2}-R+(k-1)D) <0,$$
and so the claim follows from the exact sequence. To see this, recall that
$$\deg K_{C_2}=2g_{C_2}-2=d_2^2-3d_2,$$
where $g_{C_2}$ is the genus of the smooth curve $C_2$.
\endproof

\begin{rk}
\label{rk20} The formula for $\psi_k$ given in Theorem \ref{thm20} implies the following fact: for any $h \in \wJ_{f_1}$, one has $f_2^2h \in 
 \wJ_{f}$. When $C_1$ is a smooth curve, then $ \wJ_{f_1}=S$ and this situation occurs already in 
\cite[Proposition 3.2. (i)]{BST}, where $S$ is the  polynomial ring in $n$-variables with coefficients in an arbitrary infinite field.

\end{rk}

\begin{cor}
\label{cor20}
With the above notation, assume that $C_2$ is a smooth curve, and that all singularities of $C_1$ and $C$ are quasihomogeneous.  Let $R$ be the reduced scheme of $C_1 \cap C_2$. If
$$|R|>(r_1+1)d_2,$$
then $r=r_1+d_2$.
This applies in particular when $C_2$ is a generic  curve and $r_1\ne d_1-1$.
\end{cor} 
\proof
The hypothesis $|R|>(r_1+1)d_2$ implies that the inequality
$$(k+2)d_2 <d_2^2+ |R|$$
holds for all $k \leq r_1+d_2-1$. Using Theorem \ref{thm20} and the definition of $r_1$, it follows that $D_0(f)_k=0$ for all $k \leq r_1+d_2-1$.
The exact sequence in Theorem \ref{thm20} also implies that 
$D_0(f)_{r_1+d_2} \ne 0$, which proves our claim. When $C_2$ is a generic  curve, then $C_1 \cap C_2$ consists of $d_1d_2$ nodes for $C$ and the claim is clear.
\endproof
\begin{ex}
\label{ex20}
Let $C_1$ be a reduced curve satisfying $r_1 \ne d_1-1$ and such that all the singularities of $C_1$ are quasihomogeneous. Then for any point $p \notin C_1$  and  any  line $C_2$ through $p$ such that $C_2$ meets transversally $C_1$ at smooth points, one has $r=r_1+1$. The claim follows from Corollary \ref{cor20}, since $d_2=1$ and $C_1 \cap C_2$ consists of $d_1$ nodes for $C$ in this case.
This result should be compared to Theorem \ref{thm1}. Moreover, Example \ref{ex1} shows that the restriction $r_1 \ne d_1-1$
is necessary.
If we take $p \in C_1$, then the condition that $(C,p)$ is quasihomogeneous {\it limits drastically the choices} for the line $C_2$ passing through $p$.
Consider the rational cuspidal curve $C_1:f_1=xy^{d_1-1}+z^{d_1}=0$, with $d_1>2$. We have seen in  Example \ref{ex2} that, if we take $C_2$ to be the line through the  singular point $p=(1:0:0)$  given by 
$y=0$ or $z=0$, then $(C,p)$ is quasihomogeneous, and $r=r_1$ in these two cases. In fact, in these cases $|R|=1$ and Corollary \ref{cor20} does not apply. When $C_2$ is given by $sy+tz=0$ with $st \ne 0$, then
the singularity $(C,p)$ is not quasihomogeneous, as we have seen in 
Example \ref{ex2} for the case $s=t=1$. 
\end{ex}
Assume from now on that $|R|\leq (r_1+1)d_2$, or equivalently $(k+2)d_2 \geq d_2^2+ |R|$ and set 
$$k_0=d_2-2+ \left \lceil \frac{ |R|}{d_2}\right \rceil.$$
To simplify the discussion, we also assume that $C_2$ is a smooth rational curve, hence $d_2 \in \{1,2\}$. It follows that $k_0$ is the smallest integer $k$ such that $H^0(C_2,\OO_{C_2}(-K_{C_2}-R+(k-1)D))\ne 0$. If we assume in addition that $C_1$ is a free curve, then $N(f_1)=0$ and Theorem \ref{thm20}  implies the following.
\begin{cor}
\label{corF}
With the above notation and assumptions, if in addition $C_1$ is a free curve and $C_2$ is rational, then $|R|\leq (r_1+1)d_2$ implies $k_0\geq r_1$ and 
$$r_1 \leq r= k_0 \leq r_1+d_2-1.$$
In particular, $|R| > (r_1+1)d_2-d_2^2$, that is we have the following cases.
\begin{enumerate}

\item Let $C_1:f_1=0$ be a free curve and $L$ be a line such that $C_1$ and $C_1 \cup L$ have only quasihomogeneous singularities. Then
$$|C_1 \cap L| > r_1=mdr(f_1).$$

\item  Let $C_1:f_1=0$ be a free curve and $Q$ be a smooth conic such that $C_1$ and $C_1 \cup Q$ have only quasihomogeneous singularities. Then
$$|C_1 \cap Q| > 2r_1-2, \text{ where } r_1=mdr(f_1).$$

\end{enumerate}

\end{cor}
\proof
Note that $r \geq r_1$  implies $k_0 \geq r_1$, which yields in particular the last claim.
\endproof 

When $C_2$ is a line, then $r=r_1$ and $|R|=r_1+1$ in these conditions, a known result when $C_1$ is a line arrangement, see for instance \cite[Theorem 3.6 (2)]{ADS}.
\begin{ex}
\label{exF}
Consider $C_1: xyz(x-y)(y-z)(x-z)=0$, which is free with $d_1=6$ and $r_1=2$. Let $C_2$ be a general conic passing through 2 triple points and 2 double points of $C_1$, for instance
$C_2:x^2+z^2-xy-yz=0$. Then $d_2=2$, $r_2=1$ and $|R|=6$. Corollary \ref{corF} implies
$r=k_0=3.$ It follows that the curve $C$ is free with exponents $(3,4)$ by \eqref{F1}.
Indeed, this curve $C$ has 3 nodes, 4 ordinary triple points and 2 ordinary quadruple points, hence $\tau(C)=37$.
\end{ex}
The application to the exact sequence \eqref{e8} to study free curves goes back to \cite{STY}. Now we show that this sequence gives valuable information even when $C_1$ is not a free curve. We start with the case $C_2$ is a line, hence we have to decide by Proposition \ref{prop1} or by Theorem \ref{thm1} whether $r=r_1$ or $r=r_1+1$.
\begin{cor}
\label{cor21}
With the above notation, assume that  $C_2$ is a line  and that all singularities of $C_1$ and $C$ are quasihomogeneous.  Let $R$ be the reduced scheme of $C_1 \cap C_2$. If 
$$|R|\leq r_1+1,$$
then there is the following exact sequence 
$$ 0 \to   D_0(f)_{r_1} \to H^0(C_2,\OO_{C_2}(r_1+1-|R|)) 
 \to N(f_1)_{r_1+d_1-2} \to N(f)_{r_1+d_2+d-2}   \to 0. $$
 \end{cor} 
\proof
The proof is as above, using Theorem \ref{thm20} for $k =r_1$
and the fact that $H^1(C_2,\OO_{C_2}(\ell)) =0$ if $\ell \geq 0$.
\endproof

\begin{cor}
\label{cor21NF}
\begin{enumerate}

\item Let $C_1:f_1=0$ be a nearly free curve and $L$ be a line such that $C_1$ and $C_1 \cup L$ have only quasihomogeneous singularities. Then
$$|C_1 \cap L| \geq  r_1=mdr(f_1).$$

\item  Let $C_1:f_1=0$ be a nearly free curve and $Q$ be a smooth conic such that $C_1$ and $C_1 \cup Q$ have only quasihomogeneous singularities. Then
$$|C_1 \cap Q| \geq 2r_1-1, \text{ where } r_1=mdr(f_1).$$

\end{enumerate}

 \end{cor} 
\proof
The proof is as above, using Theorem \ref{thm20} for $k =r_1-1$
and the fact that $\dim H^(C_2,\OO_{C_2}(\ell)) \geq 2> \nu(C_1)=1$ if $\ell \geq 1$. For $d_2=2$, we use the stronger fact that $\sigma(C_1)=d_1+r_1-3$.

\endproof

\begin{ex}
\label{ex21}
Let $C_1$ be a nearly free curve having only quasihomogeneous singularities.   Let $C_2$ be a line such that
$|R| =r_1$, the minimal possible value, and $C=C_1 \cup C_2$ has again only quasihomogeneous singularities. Then in the exact sequence of Corollary \ref{cor21} we have
$$ \dim H^0(C_2,\OO_{C_2}(r_1+1-|R|))\geq 2=\dim H^0(C_2,\OO_{C_2}(1)) >\dim N(f_1)_{r_1+d_1-2}=1.$$
The last equality follows from the equality $\sigma(C_1)=d_1+r_1-3$, see \cite[Corollary 2.17]{DStRIMS}. This implies $r=r_1$ in this situation.\\
A first explicit example of such a situation is provided by the curves $C_1$ discussed in Example \ref{ex2} with the line $C_2$ given by $y=0$ or $z=0$, when $r_1=1$. \\
A second  example  is provided by
the quartic with 3 cusps
$$C_1: x^2y^2+y^2z^2+x^2z^2-2xyz(x+y+z)=0,$$
which is nearly free with  $r_1=2$, see \cite[Example 2.13]{DStRIMS} and $C_2:z=0$, a line joining 2 cusps.
The curve $C$ has in this case one cusp $A_2$ and two $D_5$ singularities, hence has  only quasihomogeneous singularities.
Since $|R|=2=r_1$, the above discussion applies and
it follows that $r=r_1=2$.
 Using \cite{Dmax,dPW}, it follows that  the obtained quintic curve 
$$C: x^2y^2z+y^2z^3+x^2z^3-2xyz^2(x+y+z)=0,$$ 
 is free with exponents $(2,2)$.\\
 As a third example, consider $C_1$ to be the union of two smooth conics tangent to each other in two points, as in Proposition \ref{prop10}.
 Let $C_2$ be the line joining these two points. Then $C$ has two $D_6$ singularities, $2=|R|>r_1=1$. Hence {\it Corollary \ref{cor21} cannot be used to conclude}. Note that using the equation \eqref{M1}, we see that 
 $y\partial_x+x\partial_y \in D(f_1)_1 \cap D(f_2)_1 \ne 0$, and hence $r=1$.
  Using \cite{Dmax} it follows that this curve $C$ is nearly free with exponents $(1,4)$.
\end{ex}

Here is the version of Corollary \ref{cor21} when $C_2$ is a smooth conic. Here we know already that $r_1 \leq r \leq r_1+2$ by Theorem \ref{thm10}.
\begin{cor}
\label{cor21C}
With the above notation, assume that  $C_2$ is a smooth conic and that all singularities of $C_1$ and $C$ are quasihomogeneous.  Let $R$ be the reduced scheme of $C_1 \cap C_2$. If 
$$|R|\leq 2(r_1+1),$$
then there is the following exact sequences. 
$$ 0 \to   D_0(f)_{r_1} \to H^0(C_2,\OO_{C_2}(2r_1-|R|)) 
 \to N(f_1)_{r_1+d_1-3} \to N(f)_{r_1+d_2+d-3}  $$
 and
 $$ 0 \to   D_0(f)_{r_1+1} \to H^0(C_2,\OO_{C_2}(2r_1+2-|R|)) 
 \to N(f_1)_{r_1+d_1-2} \to N(f)_{r_1+d_2+d-2}   \to 0.$$
 \end{cor} 
 \proof
Use Theorem \ref{thm20} for $k =r_1$ and for $k=r_1+1$.
\endproof

\begin{ex}
\label{ex10bis}
Consider next the case  $C_1:f_1=xyz=0$ and $C_2:f_2=xy+yz+xz=0$. 
Then $C_2$ is a smooth conic circumscribed in the triangle $C_1$, as in the second part of Example \ref{ex10}. In this case $r_1=1$ and $|R|=3$,
hence we can apply Corollary \ref{cor21C}. The first exact sequence implies that $D_0(f)_1=0$, and the second exact sequence implies that
$$\dim D_0(f)_2=2=\dim H^0(C_2,\OO_{C_2}(1)) $$
since $C_1$ is a free curve, and hence $N(f_1)=0$.
\end{ex}
\begin{ex}
\label{ex11bis}
Let $C_1:f_1=(x-z)(3y^2-(x+2z)^2) (x^2+y^2-4z^2)=0$ be a smooth conic $Q$ circumscribed in a triangle $\Delta$ as in Example \ref{ex10bis}.
 Let $C_2:f_2= (x^2 + y^2 -z^2)=0$ be a conic inscribed in the triangle $\Delta$ and tangent to the conic $Q$ in two points. Then $d_2=2$, $r_2=1$. In this case $r_1=2$ and $|R|=5$,
hence we can apply Corollary \ref{cor21C}. The first exact sequence implies that $D_0(f)_2=0$, and the second exact sequence implies that
$$\dim D_0(f)_3=3=\dim H^0(C_2,\OO_{C_2}(1)) $$
since $C_1$ is a free curve, and hence $N(f_1)=0$.
\end{ex}

\begin{ex}
\label{ex21C}
Let $C_1$ be a nearly free curve having only quasihomogeneous singularities.   Let $C_2$ be a smooth conic such that either $|R| \leq 2r_1-1$ or
$|R| =2r_1+1$ and $C=C_1 \cup C_2$ has again only quasihomogeneous singularities. Then, when  $|R| = 2r_1-1$, we get
exactly as in Example \ref{ex21} that $r=r_1$. Assume now that $|R| =2r_1+1$.
In the first exact sequence of Corollary \ref{cor21C} we have 
$$H^0(C_2,\OO_{C_2}(2r_1-|R|))=H^0(C_2,\OO_{C_2}(-1))=0,$$
and hence $ D_0(f)_{r_1} =0$.
 In the second exact sequence of Corollary \ref{cor21C} we have
$$2 = \dim H^0(C_2,\OO_{C_2}(2r_1+2-|R|)) =H^0(C_2,\OO_{C_2}(1)) >\dim N(f_1)_{r_1+d_1-3}=1.$$
The last equality follows from the equality $\sigma(C_1)=d_1+r_1-3$, see \cite[Corollary 2.17]{DStRIMS}. This implies $r=r_1+1$ in this situation.\\
To have an explicit example, we consider again the quartic $C_1$ with 3 cusps from Example \ref{ex21}, and take now $C_2$ to be a smooth generic conic passing through the 3 cusps, then the resulting curve $C$ will have 3 $D_5$ singularities and 2 nodes $A_1$.
 It follows that $|R|=5=2r_1+1$. It follows that in this case
 $r=r_1+1=3$. As an explicit example, one can take
$$C: (xy+yz+xz)(x^2y^2+y^2z^2+x^2z^2-2xyz(x+y+z))=0,$$ 
which is a plus one generated curve with exponents $(3,3,4)$. 
\end{ex}
\begin{ex}
\label{ex21D}
Let $C_1$ be a nearly free curve having only quasihomogeneous singularities.   Let $C_2$ be a smooth conic such that  $|R|= 2r_1+2$.  In the first exact sequence of Corollary \ref{cor21C} we have 
$$\dim H^0(C_2,\OO_{C_2}(2r_1-|R|))=\dim H^0(C_2,\OO_{C_2}(-2))=0,$$
and in the second 
exact sequence of Corollary \ref{cor21C} we have 
$$\dim H^0(C_2,\OO_{C_2}(2r_1+2-|R|))=\dim H^0(C_2,\OO_{C_2})=1.$$
If  $N(f_1)_{r_1+d_1-2}=0$, then we have $r=r_1+1$.
We consider the following explicit situation. Let $C_1: f_1=x^2y^2+z^4-xz^3-2xyz^3=0$, which has an $A_4$-singularity at $p=(0:1:0)$ and an $A_2$-singularity at $q=(1:0:0)$. Then $C_1$ is a nearly free curve with 
$d_1=4$ and $r_1=2$, see \cite[Example 2.13]{DStRIMS}. Let $C_2: xy+yz+xz=0$ be a smooth conic, passing transversely through $p$ and $q$ and meeting $C_1$ transversally at another 4 points. It follows that
$|R|=2+4=6=2r_1+2$ and $C$ has 4 nodes $A_1$, one $D_5$ singularity and one $D_7$-singularity in all. Note that 
$$\sigma(C_1)=d_1+r_1-3=T_{f_1}/2=3.$$
This implies that $N(f_1)_{4}=0$, exactly what we need to conclude that $r=r_1+1=3$.
\end{ex}

\begin{ex}
\label{ex21E}
We end with an example in order to conclude via Corollary \ref{cor21C} we need to analyze the morphism  
$$N(f_1)_{r_1+d_1-3} \stackrel{f_2^2} \longrightarrow N(f)_{r_1+d-1}.$$
Let $C_1: f_1=(x^2-2xz+y^2)(x^2+2xz+y^2)=0$, hence $C_1$ is a pair of smooth conics tangent at one point, as in Proposition \ref{prop10} (2). 
This curve $C_1$ is a plus one generated curve with exponents $(2,2,3)$. Hence $d_1=4$ and $r_1=2$. Let $C_2: x^2+y^2-4z^2=0$, a circle tangent to each circle in $C_1$ at one point and passing through the 2 nodes of $C_1$. Hence $C$ has 3 singularities $A_3$ and 2 singularities $D_4$. It follows that $|R|=4=2r_1$. In the first exact sequence of Corollary \ref{cor21C} we have 
$$\dim H^0(C_2,\OO_{C_2}(2r_1-|R|))=\dim H^0(C_2,\OO_{C_2})=1,$$
and $\dim N(f_1)_{r_1+d_1-3}=2$. More precisely, a basis of $N(f_1)_{r_1+d_1-3}$ is given by the monomials $xyz$ and $xz^2$.
A computation using Singular shows that the kernel of the morphism
$$N(f_1)_{r_1+d_1-3} \stackrel{f_2^2} \longrightarrow N(f)_{r_1+d-1}$$
is 1-dimensional, generated by $xyz-xz^2$. As a result $D_0(f)_2=0$.
 In the second 
exact sequence of Corollary \ref{cor21C} we have 
$$\dim H^0(C_2,\OO_{C_2}(2r_1+2-|R|))=\dim H^0(C_2,\OO_{C_2}(2))=3$$
and $\dim N(f_1)_{r_1+d_1-2}=1$. Therefore we have $r=r_1+1=3$.
A computation with Singular shows that $C$ is a plus one generated curve with exponents $(3,3,4)$.
\end{ex}

\subsection{An application to the jumping lines of the rank 2 vector bundle $E_{C_1}$}

For a reduced plane curve $C$ and a line $L$ in $\PP^2$, the pair of integers $( d_1^L(C), d_2^L(C))$ such that  $ d_1^L(C) \leq d_2^L(C)$  and  $ E_C|_L \simeq \OO_L(-d_1^L(C)) \oplus \OO_L(-d_2^L(C))$ is called the (ordered) splitting type  of $E_C$ along $L$, see for instance \cite{OSS}. 
For a generic line $L_0$, the corresponding splitting type $( d_1^{L_0}(C), d_2^{L_0}(C))$ is known to be constant, see \cite[Definition 2.2.3 and Lemma 3.2.2]{OSS}.
A line $L$ in $\PP^2$ is called a jumping line  for $E_C$ or, equivalently, for $T\langle C \rangle $, if 
$$d_1^{L_0}(C)-d_1^{L}(C)>0.$$
The following result relates the splitting type of $E_C$ along a line $L: \al_L=0$, to the Lefschetz properties of the Jacobian module $N(f)$ with respect to the multiplication by $\al_L$, see \cite[Proposition 4.1]{DStJump}.
\begin{prop} \label{propthmA} 
 For any reduced curve $C:f=0$ and any line $L: \al_L=0$ in $\PP^2$, we have
 $d_1^L(C)=\min \{mdr(f), k(f,L)\},$
 where $$k(f,L)= \min \{ k \in \N \ : \ N(f)_{k+d-2} \stackrel{\cdot \alpha_L}{\to}
N(f)_{k+d-1} \text{ is not injective } \}.$$
 \end{prop}
 Using this result, we give now an easy geometric way to check that a line is a jumping line, under some conditions.
\begin{thm} \label{thmJ} 
Let $C_1:f_1=0$ be a reduced curve and $C_2:f_2=0$ be a line in $\PP^2$. 
 Assume that  all the singularities of $C_1$ and of $C=C_1 \cup C_2$ are quasihomogeneous, and let $R$ be the reduced scheme of $C_1 \cap C_2$. If $|R| < r_1+1$, then the morphism
 $$ \psi'_k: N(f_1)_{k+d_1-2} \stackrel{f_2} \longrightarrow N(f_1)_{k+d_1-1} $$
 is not injective, for any $k$ satisfying $|R|-1 \leq k <r_1$ and one has
 $$d_1^{C_2}(C_1)=k(f_1,C_2) \leq |R|-1.$$
 Moreover, if one of the following two conditions holds
 \begin{enumerate}

\item  either $2r_1<d_1$, or

\item  $$ 2r_1 \geq d_1 \text{ and } |R|-1 < \left \lfloor \frac{d_1-1}{2} \right \rfloor,$$
 
\end{enumerate}
then $C_2$ is  a jumping line for  the rank two vector bundle $E_{C_1}$.

\end{thm}
\proof
Note that the condition $k \geq |R|-1$ implies that 
$$H^0(C_2,\OO_{C_2}(-K_{C_2}-R+(k-1)D))=H^0(C_2,\OO_{C_2}(k+1-|R|))\ne 0$$
in the exact sequence from Theorem \ref{thm20}. On the other hand, the condition $k<r_1$ implies that $D_0(f)_k=0$. 
Hence, using  Theorem \ref{thm20}, we see that the morphism
$$\psi_k: N(f_1)_{k+d_1-2} \stackrel{f_2^2} \longrightarrow N(f)_{k+d_1}.$$
is not injective. To prove our claim,  it is enough to show that 
$$\ker  \psi_k \subset \ker \psi'_k.$$
Let $h \in \ker \psi_k$ be (the representative of ) some element in this kernel. This means that $f_2^2h \in J_f$, in other words there is a derivation $\delta \in Der(S)$ such that
$$f_2^2h=\delta (f)= f_2 \delta (f_1)+ f_1 \delta (f_2).$$
Since $f_1$ and $f_2$ have no common factor, this means that $\delta (f_2)$ is divisible by $f_2$, say $\delta (f_2)=f_2 g$ for some $g \in S$.
Dividing the above relation by $f_2$ we get
$$f_2h=  \delta (f_1)+ f_1g,$$
which implies $f_2h \in J_{f_1}$. Hence $h \in \ker \psi'_k$ as we claimed. Now we prove the final claims in Theorem \ref{thmJ}.
Since we know that $\psi'_k$ is not injective for $|R|-1 \leq k <r_1$,
it follows by Proposition \ref{propthmA} that 
$$d_1^{C_2}(C_1)=k(f_1,C_2) \leq |R|-1.$$
Assume now $2r_1<d_1$. First note that the curve $C_1$ cannot be a free curve in view of Corollary \ref{corF} $(1)$ saying that in this case $|R|>r_1$. If $C_1$ is a nearly free curve, then the exponents $r=d_1'\leq d_2'$ verify $d_1'+d_2'=d_1$, and hence the condition $2r_1<d_1$ implies $d_1'<d_2'$. Using for instance \cite[Example 4.8]{DStJump}
we see that in this case $d_1^{L_0}(C_1)=r_1$. The same equality holds for all the other non free reduced plane curves $C_1$ satisfying $2r_1<d_1$, see \cite[Corollary 4.5]{DStJump}. This fact implies that 
$C_2$ is a jumping line for $E_{C_1}$ in the case $(1)$. The claim in case (2) follows from  \cite[Corollary 4.6 and Example 4.8]{DStJump}.
Indeed, these results imply that we have
$$d_1^{L_0}(C_1)= \left \lfloor \frac{d_1-1}{2} \right \rfloor$$
when $2r_1 \geq d_1$.
\endproof

\begin{ex}
\label{ex1.1} Let $C_1: f_1=0$ be a Thom-Sebastiani plane curve, i.e. a curve such that $f_1(x,y,z)=g(x,y)+z^{d_1}$, where
$g$ is a homogeneous polynomial of degree $d_1$ in $S'=\C[x,y]$. Assume that
$$g=\ell_1^{k_1} \cdots \ell_m^{k_m},$$
where the linear forms $\ell_j \in S'$ are distinct and $m\geq 2$. 
It follows that $C_1$ is a 3-syzygy curve with exponents $r_1=d'_1=m-1$ and $d'_2=d'_3=d_1-1$ for $m \geq 3$ and $C_1$ is nearly free with exponents $r_1=d'_1=1$ and $d'_2=d_1-1$ for $m=2$, see \cite[Example 4.5]{DSt3syz}. 

Assume first that $2(m-1)<d_1$ and take the line $C_2$ to be given by one of the factors of $g$, say $C_2:f_2=\ell_1=0$. Then it is easy to check that all the singularities of $C_1$ and of $C=C_1 \cup C_2$ are quasihomogeneous. To do this, one can assume that $\ell_1=x$ and hence $R$ is the point $(0:1:0)$. It follows that $|R|=1<r_1+1$, and Theorem \ref{thmJ} (1) implies that the line $C_2$ is a jumping line for $E_{C_1}$ and moreover
$$d_1^{C_2}(C_1)=k(f_1,C_2)=0.$$
Note that the inequality $d_1^{C_2}(C_1)=k(f_1,C_2) \geq 0$ holds in general, see for instance \cite[Proposition 2.5]{DStJump}. 

If we assume now that $2(m-1)\geq d_1\geq 3$, then we get the same result using Theorem \ref{thmJ} (2).
In this way we have found out $m$ jumping lines $\ell_j=0$ for $j=1, \ldots, m$ for the vector bundle $E_{C_1}$.

\end{ex}
Recall that when $C_1$ is a free curve, then all the lines $L$ in $\PP^2$ are not jumping lines for the vector bundle $E_{C_1}$. With this in mind, the following result, which is a reformulation of Theorem \ref{thmJ}, can be regarded as a generalization of Corollary \ref{corF} (1).
\begin{cor}
\label{corR} 
Let $C_1:f_1=0$ be a reduced curve and $L$ be a line in $\PP^2$, which is not a jumping line for the vector bundle $E_{C_1}$. Let $d_1=\deg(f_1)$ and $r_1=mdr(f_1)$. If all the singularities of $C_1$ and of $C_1 \cup L$ are quasihomogeneous, then
$$|C_1 \cap L| >d_1^{L}(C_1)= \begin{cases} r_1 &\mbox{if }  \  2r_1<d_1, \\ 
 \left \lfloor \frac{d_1-1}{2} \right \rfloor & \mbox{if }  \  2r_1 \geq d_1. \end{cases} $$
\end{cor}

\end{document}